\documentclass{article}
\usepackage{amssymb,amsthm}
\usepackage{graphicx} 
\usepackage{xcolor}

\newtheorem{proposition}{Proposition}

\title{A semifield of order 128 and fractional dimension $\frac{7}{3}$ relative to one of its subsemifields}
\author{I.F. R\'ua\footnote{Department of Mathematics, University of Oviedo, Asturias, Spain. {rua@uniovi.es}}\ ,\ E.F. Combarro\footnote{Department of Computer Science, University of Oviedo, Asturias, Spain.}}
\date{April 2, 2025}

\begin{document}

\maketitle

\begin{abstract}
    In this short note, an example of a semifield of order 128 containing the Galois field $\mathbb{F}_8$ is given. Up to our knowledge, this is the first example supporting the following problem by Cordero and Chen (2013): ``There exist semifield planes of order $2^t$, for any integers $t$ relatively prime to 3 that admit semifield subplanes of order $2^3$''. 
\end{abstract}
\section{Introduction}

A finite semifield $(S,+,\cdot)$ can be seen as a finite nonassociative division algebra, or as a special kind of ternary ring \cite{Knu65,CorWene,Lavrauw}. From this second point of view, it coordinatizes a projective plane or order $|S|=p^n$, for some prime number $p$. When $S$ contains a subsemifield $T$ of order $p^m$ (i.e., a subset which is a semifield under the restriction of $+$ and $\cdot$), then the dimension of $S$ with respect to $T$ is $\frac{n}{m}$. When $S$ is a finite field, it is well-known that such a dimension can only be an integer number. However, for proper semifields (i.e., non-associative), which coordinatized non-Desarguesian planes, it can be a non-integer rational number. In such a case, it is called a fractional (dimensional) semifield \cite{Cordero1,Cordero2,Olga,Trombetti}.      

Up to our knowledge, any fractional semifield known so far has order $|S|=p^n$ (with $p=2,3$), with the corresponding subsemifield of order $|T|=p^2$ \cite{Jha}. In \cite{Chen}, it was posed the following problem: ``There exist semifield planes of order $2^t$, for any integers $t$ relatively prime to 3 that admit semifield subplanes of order $2^3$'', claiming that ``no examples to support it have been found''.

In this short note, an example of a semifield of order 128 containing the Galois field $\mathbb{F}_8$ is given. So, it is a fractional semifield of dimension $\frac{7}{3}$. This semifield was found during 
an ongoing investigation of the structure of semifields of order 128.

\section{The example} 
A finite semifield $S$ (i.e., a finite nonassociative division ring) with $2^7$ elements is a $7-$dimensional algebra over the Galois field $\mathbb{F}_2$ with unity $1$. If $\mathcal B=\{x_1,\dots,x_7\}$ is a $\mathbb{F}_2$-basis of $S$, then there exists a unique set of constants (known as {cubical array} or {$3-$cube} corresponding to $S$ with respect to the basis $\mathcal B$) $\mathbf A_{S,\mathcal B}=\{A_{i_1i_2i_3}\}_{i_1,i_2,i_3=1}^7\subseteq \mathbb{F}_2$ such that 
$$x_{i_1}x_{i_2}=\sum_{i_3=1}^7{A_{i_1i_2i_3}}x_{i_3}\; \forall i_1,i_2\in\{1,\dots,7\},$$
that fully determines the multiplication in $S$. This fact allows a concrete description of semifields of order 128 by specific sets of $7\times 7$ matrices over the field of $2$ elements.

\begin{proposition}\label{matrices}\cite{Hentzel}[Proposition 3]
    For any finite semifield $S$ of order $2^7$, there exists a
    set of $7$ matrices, called \emph{standard basis}, $\{A_i\}_{i=1}^7\subseteq GL(7,2)$, such 
    that:
    \begin{enumerate}
        \item $A_1$ is the identity matrix;
        \item Any non-zero $(\lambda_1,\dots,\lambda_7)$, the linear combination $\sum_{i=1}^7\lambda_i A_i$ is invertible;
        \item The first column of the matrix $A_i$ is the column vector
        $e_i^\downarrow$ (with a $1$ in
    the $i$-th position, and $0$ everywhere else).
    \end{enumerate}
    Such a standard basis is the set of coordinate matrix of the $\mathbb{F}_2-$linear maps $\{L_{a_i}(w)=a_iw\}_{i=1}^7$, for a $\mathbb{F}_2-$basis $\mathcal B=\{a_1=1,a_2,\dots,a_7\}$ of $S$. In terms of the $3-$cube $\mathbf A_{S,\mathcal B},$ we have $(A_i)_{kj}=A_{ijk},$ for all $i,j,k=1,\dots,7$.
\end{proposition}

Next, we provide a standard basis for a semifield $S$ of order 128 containing the Galois field $\mathbb{F}_8$ as a subsemifield. Namely,

$$A_1=\left(\begin{array}{ccc|cccc}1&0&0&0&0&0&0\\0&1&0&0&0&0&0\\0&0&1&0&0&0&0\\\hline0&0&0&1&0&0&0\\0&0&0&0&1&0&0\\0&0&0&0&0&1&0\\0&0&0&0&0&0&1\end{array}\right)$$
$$A_2=\left(\begin{array}{ccc|cccc}0&0&1&0&0&0&0\\1&0&1&0&0&0&0\\0&1&0&0&0&0&0\\\hline0&0&0&1&1&1&1\\0&0&0&1&0&1&0\\0&0&0&0&0&1&1\\0&0&0&0&0&1&0\end{array}\right)\ , \ A_3=\left(\begin{array}{ccc|cccc}0&1&0&0&0&0&0\\0&1&1&0&0&0&0\\1&0&1&0&0&0&0\\\hline0&0&0&0&1&0&0\\0&0&0&1&1&0&0\\0&0&0&1&1&1&1\\0&0&0&1&0&1&0\end{array}\right)$$$$A_4=
\left(\begin{array}{ccccccc}0&0&0&1&1&0&0\\0&0&0&1&1&1&0\\0&0&0&0&1&0&1\\1&0&0&0&0&0&1\\0&0&1&0&0&1&0\\0&1&1&1&0&1&1\\0&0&0&1&1&0&1\end{array}\right)\ , \ A_5=\left(\begin{array}{ccccccc}0&0&0&0&1&0&0\\0&0&0&0&1&1&1\\0&0&0&1&0&1&0\\0&0&1&0&0&1&1\\1&0&1&0&0&0&1\\0&0&0&0&1&1&0\\0&1&1&1&0&1&1\end{array}\right)$$$$A_6=
\left(\begin{array}{ccccccc}0&0&0&1&0&1&1\\0&0&0&1&1&0&0\\0&0&0&0&1&0&0\\0&0&0&1&0&1&0\\0&1&1&1&1&1&1\\1&0&1&0&1&1&0\\0&0&1&1&0&1&1\end{array}\right)\ , \ A_7=
\left(\begin{array}{ccccccc}0&0&0&1&1&0&1\\0&0&0&0&1&0&0\\0&0&0&1&0&0&0\\0&1&1&0&1&0&1\\0&1&1&1&0&1&0\\0&0&1&1&1&0&1\\1&0&0&0&1&1&0\end{array}\right)$$

The fact that this set of matrices is indeed a standard basis can be straightforwardly checked with the help of Proposition \ref{matrices}. The existence of a subsemifield of order 8 can be clearly deduced from the form of the first three matrices. They are diagonal block matrices in which the upper left part corresponds to the first three powers of the companion matrix $C(x^3+x+1)$. Since that polynomial is irreducible over $\mathbb{F}_2$, the subsemifield given by the linear combinations of the elements $\{a_1=1,a_2,a_3=a_2^2\}$ (with $a_2^3=a_2+1$) is $\mathbb{F}_8$. Hence, the fractional dimension $\frac{3}{7}$.

As an aside, let us comment that this semifield $S$ is not commutative, so its opposite ring $S^{opp}$ (i.e., the one for which the standard basis is related to the $\mathbb{F}_2-$linear maps $\{R_{a_i}(w)=wa_i\}_{i=1}^7$), is also fractional. It corresponds to a semifield in the Knuth orbit of $S$, because its $3-$cube is 
$\mathbf A_{S^{opp},\mathcal B}=\{A_{i_2i_1i_3}\}_{i_1,i_2,i_3=1}^7\subseteq \mathbb{F}_2$.

\section{Conclusions}

We have given an example of a semifield of order 128 containing a subsemifield of order 8. This is the first example of a fractional semifield of dimension $\frac{n}{m},$ where $m$ is not equal to 2. As a future work, we intend to transform this example into a family of examples of fractional semifields of order $2^n$, containing $\mathbb{F}_8$ ($3\ \not|\ n$).

\section*{Acknowledgments}

This work has been partially supported by grant PID2023-146520OB-C22, funded by MCI- U/AEI/10.13039/501100011033/FEDER, EU and PID2023-146520OB-C22 funded by MICIU/AEI/10.13039/501100011033, and by the Ministry for Digital Transformation and of Civil Service of the Spanish Government through the QUANTUMENIA project call Quantum Spain project, and by the European Union through the Recovery, Transformation and Resilience Plan NextGenerationEU within the framework of the Digital Spain 2026 Agenda.
This work has been financially supported by grant PID2021-123461NB-C22, and by the Ministry of Economic Affairs and Digital Transformation of the Spanish Government through the Spanish National Institute of Cibersecurity (INCIBE) project call for Strategic Projects on Cibersecurity in Spain (Grant MRR-MAETD-24-INCIBE-001), and by the European Union through the Recovery, Transformation and Resilience Plan NextGenerationEU within the framework of the Digital Spain 2026 Agenda. The first author is a member of the Universidad de Oviedo research team GACYC, and have also been supported by the Spanish Network of Mathematics in the Information Society (MatSI).

\bibliographystyle{plain}
\bibliography{quantum}

\begin{thebibliography}{10}

\bibitem{Chen}
Linlin Chen and Minerva Cordero.
\newblock Fractional dimensional semifield planes.
\newblock {\em Note Mat.}, 32(2):57--61, 2012.

\bibitem{CorWene}
M.~Cordero and G.~P. Wene.
\newblock A survey of finite semifields.
\newblock volume 208/209, pages 125--137. 1999.
\newblock Combinatorics (Assisi, 1996).

\bibitem{Cordero1}
Minerva Cordero and Vikram Jha.
\newblock Primitive semifields and fractional planes of order {$q^5$}.
\newblock {\em Rend. Mat. Appl. (7)}, 30(1):1--21, 2010.

\bibitem{Cordero2}
Minerva Cordero and Vikram Jha.
\newblock Fractional dimensions in semifields of odd order.
\newblock {\em Des. Codes Cryptogr.}, 61(2):197--221, 2011.

\bibitem{Hentzel}
I.~R. Hentzel and I.~F. R\'ua.
\newblock Primitivity of finite semifields with 64 and 81 elements.
\newblock {\em Internat. J. Algebra Comput.}, 17(7):1411--1429, 2007.

\bibitem{Jha}
Vikram Jha.
\newblock The ubiquity of the orders of fractional semifields of even
  characteristic.
\newblock {\em Des. Codes Cryptogr.}, 72(3):675--686, 2014.

\bibitem{Trombetti}
Vikram Jha, Olga Polverino, and Rocco Trombetti.
\newblock Subplanes of a translation plane.
\newblock {\em Finite Fields Appl.}, 30:121--138, 2014.

\bibitem{Knu65}
D.~E. Knuth.
\newblock Finite semifields and projective planes.
\newblock {\em {\rm Journal of Algebra}}, 2:182--217, 1965.

\bibitem{Lavrauw}
Michel Lavrauw.
\newblock Finite semifields and nonsingular tensors.
\newblock {\em Des. Codes Cryptogr.}, 68(1-3):205--227, 2013.

\bibitem{Olga}
Olga Polverino and Rocco Trombetti.
\newblock Fractional dimension of binary {K}nuth semifield planes.
\newblock {\em J. Combin. Des.}, 20(7):317--327, 2012.

\end{thebibliography}

\end{document}